\documentclass[a4paper]{article} 
\usepackage[english]{babel}
\usepackage[leqno,fleqn]{amsmath}
\usepackage{amsfonts,amssymb,amsthm}
\usepackage[latin1]{inputenc} 
\usepackage[T1]{fontenc}

\newtheorem{lemma}{Lemma}
\newtheorem{satz}[lemma]{Theorem}
 
\newcommand{\nc}{\newcommand}
\nc{\Ind}{
 \setbox0=\hbox{$x$}\kern\wd0\hbox to 0pt{\hss$
 \mid$\hss}\lower.9\ht0\hbox to 0pt{\hss$\smile$\hss}\kern\wd0
}
\nc{\indep}[3]{
 #1\mathop{\mathpalette\Ind{}}_{#2}#3
}
\nc{\mr}{\mathrm{MR}}
\nc{\acl}{\mathrm{acl^{eq}(0)}}
\nc{\s}{\mathrm{Stab}}
\renewcommand{\hom}{\mathrm{Hom}}

\title{A note on generic types}
\author{Martin Ziegler} \date{August 16, 2006}
\begin{document}
\maketitle 
\begin{abstract}
  \noindent In a stable abelian group, we characterize generic types of
  cosets of type-definable subgroups.
\end{abstract}
The following remark is part of the folklore. It
was stated by the author 1990 in a letter to Daniel Lascar.
\begin{satz}
  Let $G$ be a stable abelian group. $a$, $b$ and $c$
  pairwise independent over $0$ and $a+b+c=0$. Then:
  \begin{enumerate}
  \item The strong types of $a$, $b$ and $c$ all have the same
    connected stabilizer $U$.
  \item $a$, $b$, and $c$ are generic elements of\/
    $\acl$--definable cosets of $U$.
  \end{enumerate}
\end{satz}
If $G$ is totally--transcendent, it follows that $a$, $b$ and $c$ have
the same Morley rank over $0$, namely den rank of $U$. Moreover $U$ is
definable.\\

Let $p$ and $q$ be strong types over $0$. Then
$$\hom(p,q)=\{g\in G \mid \forall a\;\; a\models p_{|g} \Rightarrow
a+g\models q_{|g}\} $$
is an $\acl$--type definable (if $G$ is totally
transcendent: definable) subset of $G$.  Of course
$\hom(p,p)=\s(p)$.
\begin{lemma}Let $p$, $q$ and $r$ be  strong types. Then
  \begin{itemize}
  \item For all $A\subset G$ and $g\in\hom(p,q)$ we have
    $$a\models p_{|g,A} \Rightarrow a+g\models q_{|g,A}.$$
  \item $\hom(p,q)+\hom(q,r)\subset\hom(p,r)$
  \item $\hom(p,q)=-\hom(q,p)=\hom(-q,-p)$
  \item $0\in\hom(p,p)$
  \item If $\hom(p,q)$ is non--empty and $G$ totally--transcendent,
    $$\mr(\hom(p,q))\leq\mr(p)=\mr(q).$$
  \end{itemize}
\end{lemma}
\begin{proof}
  Let $a$ be a realization of $p$, which is independent of $g,A$. 
  Then $\indep{a}{g}{A}$ implies that $\indep{a+g}{g}{A}$. Since
  $\indep{a+g}{}{g}$, we have $\indep{a+g}{}{g,A}$.

  Let $g\in\hom(p,q)$, $h\in\hom(q,r)$ and $a$ be a realization of $p$
  which is independent of $g+h$. We may assume that $a$ is
  independent of $g,h$.  Then $a+g$ is a realization of $q$, which
  is also independent of $g,h$. Therefore $a+g+h$ is a realization of
  $r$, which is independent of $g,h$ and therefore of $g+h$.

  Let $g\in\hom(p,q)$ and $a$ be a realization of $p$, which is
 independent of $g$.
 Then is
  $$\mr(a)=\mr(a/g)=\mr(a+g/g)=\mr(a+g).$$
  This implies $\mr(p)=\mr(q)$. From
  $$\mr(g)=\mr(g/a)=\mr(a+g/a)\leq\mr(a+g)$$
  follows $\mr(\hom(p,q))\leq\mr(q)$.
\end{proof}
Assume now that $a$, $b$ and $c$ are as in the theorem, $p$, $q$ and
$r$ the strong types of $a$, $b$ and $c$. Then, trivially,
\begin{itemize}
\item $p(G)\subset\hom(q,-r)=\hom(r,-q)$
\item $q(G)\subset\hom(r,-p))=\hom(p,-r)$
\item $r(G)\subset\hom(p,-q)=\hom(q,-p)$.
\end{itemize}
It follows
\begin{itemize}
\item $p(G)-p(G)\subset\s(q),\s(r)$
\item $q(G)-q(G)\subset\s(r),\s(p)$
\item $r(G)-r(G)\subset\s(p),\s(q)$.
\end{itemize}
Since on the other hand
\begin{itemize}
\item $\s(p)\subset p(G)-p(G)$
\item $\s(q)\subset q(G)-q(G)$
\item $\s(r)\subset r(G)-r(G)$,
\end{itemize}
all stabilizers are equal to
$$U=p(G)-p(G)=q(G)-q(G)=r(G)-r(G).$$
Therefore $p(G)$, $q(G)$, $r(G)$ all lie in
$\acl$--definable cosets of $U$. Since $U$ is the stabilizer of
each of these types, $p$, $q$, $r$ are generic types and $U$ is
connected.
This proves the theorem.
\end{document}